\documentclass[12pt]{amsart}%{article}
\usepackage{amsmath, amsthm, amssymb}
\usepackage[all]{xy}
\theoremstyle{plain}
\newtheorem{theorem}{Theorem}[section]
\newtheorem{corollary}[theorem]{Corollary}
\newtheorem{lemma}[theorem]{Lemma}
\newtheorem{proposition}[theorem]{Proposition}

\theoremstyle{definition}

\newtheorem{remark}[theorem]{Remark}

\theoremstyle{remark}

\DeclareMathOperator{\op}{op}
\DeclareMathOperator{\Ob}{Ob}

\begin{document}

\title[Presheaves of simplicial groupoids and 2-groupoids]{Closed model categories for presheaves of simplicial groupoids
and presheaves of 2-groupoids }

\author{Zhi-Ming Luo} \address{Department of Mathematics and
  Statistics, University of Guelph, Guelph, Ontario N1G 2W1 Canada}
\email{zluo@uoguelph.ca} \author{Peter Bubenik} \address{Department of
  Mathematics, Cleveland State University, 2121 Euclid Ave. RT 1515,
  Cleveland OH 44115-2214, USA} \email{p.bubenik@csuohio.edu}
\author{Peter T. Kim} \address{Department of Mathematics and
  Statistics, University of Guelph, Guelph, Ontario N1G 2W1 Canada}
\email{pkim@uoguelph.ca}



\begin{abstract}
  It is shown that the category of presheaves of simplicial groupoids and the   category of presheaves of 2-groupoids have Quillen closed model   structures. Furthermore, their homotopy categories are equivalent to the homotopy categories of simplicial presheaves and homotopy 2-types, respectively. 
\end{abstract}

\maketitle

%%%%%%%%%%%%%%%%%%%%%%%%%%%%%%%%%%%%%%%

\section{Introduction}
%%

%\bigskip

Quillen's axioms for a closed model categories have proved to be a successful framework for expanding the scope of the tools of homotopy theory. Many categories have been shown to have a useful closed model structure. A central example is
% foundational theorem of simplicial homotopy theory asserts that 
the category {\bf S} of simplicial sets~\cite{quillen:homotopicalAlgebra}. Further examples %of useful closed model structures 
related to the present study include
% Mathematicians have found a large quantity of categories enjoying the closed % model structures. For
% example, 
the category of simplicial groupoids by Dwyer and Kan
\cite{dwyerKan:simplicialGroupoids,goerssJardine:book}, the category of 2-groupoids by
Moerdijk and Svensson \cite{moerdijkSvensson:homotopy2types}, the category of simplicial
presheaves by Jardine \cite{jardine:simplicialPresheaves}, the category of simplicial
sheaves by Joyal~\cite{joyal:letter} and the category of sheaves of simplicial groupoids by Joyal and Tierney~\cite{joyalTierney:strongStacks, joyalTierney:simplicialGroupoids} and Crans~\cite{crans:sheaves},
% Crans \cite{crans:sheaves} uses adjoint
% functors to prove that a kind of sheaves have closed model
% structures according to a well-known closed model category.
and sheaves of 2--groupoids by Crans~\cite{crans:sheaves}.

We use techniques developed by Jardine~\cite{jardine:simplicialPresheaves} to prove that the
category of presheaves of simplicial groupoids and the category of
presheaves of 2-groupoids have closed model structures.

Let $\mathcal{C}$ be a small site. The adjunction $G: {\bf S} \rightleftarrows s{\bf Gd}: \bar{W}$ between the loop groupoid functor and the universal cocycle functor induces an adjunction 
\begin{equation} \label{eq:adjn1}
  G: {\bf S}Pre(\mathcal{C}) \rightleftarrows s{\bf Gd}Pre(\mathcal{C}): \bar{W}.
\end{equation}
A map $f:X \to Y \in s{\bf Gd}Pre(\mathcal{C})$ is defined to be a
fibration (weak equivalence) if $\bar{W}(f)$ is a global fibration
(topological weak equivalence) in ${\bf S}Pre(\mathcal{C})$. A map is
a cofibration if it has the left lifting property with respect to all
maps that are both fibrations and weak equivalences. We show that with
these classes of morphisms, $s{\bf Gd}Pre(\mathcal{C})$ is a closed model
category.

Similarly, there is an adjunction $G: {\bf S} \rightleftarrows {\bf 2-Gpd}: \bar{W}$, which induces an adjunction
\begin{equation} \label{eq:adjn2}
  G: {\bf S}Pre(\mathcal{C}) \rightleftarrows {\bf 2-Gpd}Pre(\mathcal{C}): \bar{W}.
\end{equation}
As with the category of presheaves of simplicial groupoids above, this
adjunction and the global fibrations and topological weak equivalences
in ${\bf S}Pre(\mathcal{C})$ induce a closed model structure on the
category of presheaves of 2-groupoids.

The model structures of Joyal and Tierney~\cite{joyalTierney:simplicialGroupoids} and Crans~\cite{crans:sheaves} are similarly obtained from a model structure on simplicial sheaves on $\mathcal{C}$, but using different functors.

We also show that the homotopy category associated to the first
category is equivalent to the homotopy category of simplicial
presheaves, and that the homotopy category associated to the latter
category is equivalent to the homotopy category of homotopy 2-types.
In fact, the adjunction \eqref{eq:adjn1} is a Quillen equivalence and
the adjunction \eqref{eq:adjn2} is a Quillen adjunction which induces
the desired equivalence of homotopy categories. 
Applications of this theory include homotopy classification results for $G$-torsors and gerbes~\cite{jardineLuo, jardine:cocycle, jardine:gerbes}

%%%%%%%%%%%%%%%%%%%%%%
%  section 2         %
%%%%%%%%%%%%%%%%%%%%%%

\section{Preliminaries}

A \emph{groupoid} is a small category in which every morphism is invertible. Let ${\bf Gd}$ denote the category of groupoids. A \emph{simplicial object} in a category ${\bf C}$ is a functor ${\bf \Delta}{^{\op}} \to {\bf C}$, where ${\bf \Delta}$ is the category of finite ordinal numbers and order-preserving maps.  A \emph{simplicial groupoid} is a simplicial object in the category of groupoids that is levelwise constant. That is, a simplicial groupoid $G$ has $\Ob(G_n) = \Ob(G)$ and for each ordinal number map $\phi: {\bf n} \to {\bf m}$, the induced map $\phi^*: G_m \to G_n$ is the identity map.  For each $x,y \in G$, the morphism sets $G_n(x,y)$ give the $n$-simplices of the simplicial set $G(x,y)$.  Let $s{\bf Gd}$ denote the category of simplicial groupoids.

We remark that our simplicial groupoids are groupoids enriched in simplicial sets. The simplicial groupoids of Joyal and Tierney~\cite{joyalTierney:simplicialGroupoids} are more general. They are the groupoids in simplicial sets.

A \emph{2-groupoid} is a strict 2-category in which every 1-morphism
and every 2-morphism has a (strict) inverse. Let ${\bf 2-Gpd}$ be the
category of small 2-groupoids and strict functors.

Given categories ${\bf C}$ and ${\bf D}$, a ${\bf D}$--valued presheaf on ${\bf C}$ is a functor ${\bf C}^{\op} \to {\bf D}$.
Morphisms of presheaves are natural transformations.
Presheaves of simplicial groupoids are $s{\bf Gd}$--valued presheaves and presheaves of 2-groupoids are ${\bf 2-Gpd}$--valued presheaves.
The category of ${\bf D}$--valued presheaves on ${\bf C}$ is denoted ${\bf D}Pre({\bf C})$. 
%A \emph{Grothendieck site} is a category with a Grothendieck topology. 

A \emph{Quillen closed model category $\mathcal{D}$} is a category
  which is equipped with three classes of morphisms, called
  cofibrations, fibrations and weak equivalences which together
  satisfy the following axioms \cite{quillen:homotopicalAlgebra}, \cite{quillen:rht}, \cite{goerssJardine:book}. Fibrations (cofibrations) that are also weak equivalences are called trivial fibrations (trivial cofibrations).
  \begin{itemize}
     \item[{\bf CM1:}] The category $\mathcal{D}$ is closed under all finite
            limits and colimits.
     \item[{\bf CM2:}] Suppose that the following diagram
            commutes in $\mathcal{D}$:
         $$\xymatrix{ X \ar[rr]^-g \ar[dr]_-h && Y\ar[dl]^-f\\
                        &   Z }
          $$
        If any two of $f,g$ and $h$ are weak equivalences, then so is the
        third.
     \item[{\bf CM3:}] If $f$ is a retract of $g$ and $g$ is a weak
         equivalence, fibration or cofibration, then so is $f$.
      \item[{\bf CM4:}] Suppose that we are given a commutative diagram
         $$\xymatrix{ U\ar[r]\ar[d]_i   & X \ar[d]^p\\
             V \ar[r]\ar@{.>}[ur]    & Y  }
            $$
        where $i$ is a cofibration and $p$ is a fibration. Then the
        lifting exists, making the diagram commute, if either $i$ or $p$
        is also a weak equivalence.
     \item[{\bf CM5:}] Any map $f: X \to Y $ may be factored:
        \begin{itemize}
           \item[(a)] $f = p \cdot i $ where $p$ is a fibration and $i$ is a trivial
                  cofibration, and
            \item[(b)] $f = q \cdot j $ where $q$ is a trivial fibration and $j$ is a
              cofibration.
         \end{itemize}
       \end{itemize}

A \emph{right proper closed model category $\mathcal{D}$} is a
closed model category satisfying the following axiom:
\begin{itemize}
\item[{\bf P1:}] the class of weak
equivalences is closed under base change by fibrations. In other
words, given a pullback diagram
$$\xymatrix{
X \ar[r]^-{g_\ast} \ar[d] & Y\ar[d]^-p\\
Z \ar[r]_-g& W}
$$
of $\mathcal{D}$ with $p$ a fibration, if $g$ is a weak
equivalence then so is $g_\ast$.
\end{itemize}

% A category $ \mathcal{D}$ is a \emph{simplicial category} if there
% is a mapping space functor
% $${\bf Hom}_{ \mathcal{D}}(\cdot, \cdot):
% \mathcal{D}^{op} \times \mathcal{D } \to \bf{S}$$ with the
% properties that for $A$ and $B$ objects in $ \mathcal{D}$
%   \begin{itemize}
%    \item[(1)] ${\bf Hom}_\mathcal{D}(A,B)_0 = \mbox{hom}_
%    \mathcal{D}(A,B)$;
%    \item[(2)] the functor ${\bf Hom}_{ \mathcal{D}}(A, \cdot):
%     \mathcal{D } \to \bf{S}$ has a left adjoint
%     $$A \otimes \cdot: {\bf S} \to \mathcal{D}$$
%     which is associative in the sense that there is an isomorphism
%     $$A\otimes(K \times L) \cong (A\otimes K) \otimes L,$$
%     natural in $A \in \mathcal{D}$ and $K,L \in {\bf S}$;
%     \item[(3)] The functor ${\bf Hom}_{ \mathcal{D}}(\cdot, B):
%     \mathcal{D}^{op}  \to \bf{S}$ has a left adjoint
%     $${\bf hom}_{ \mathcal{D}}(\cdot, B): {\bf S} \to
%     \mathcal{D}^{op}.  $$
%    \end{itemize}
% We remark that ${\bf hom}_{\mathcal{D}}(K,B)$ is often denoted $B^K$.

% A \emph{simplicial model category} $ \mathcal{D}$ is both a closed
% model category and a simplicial category which satisfies the
% following axiom:
%   \begin{itemize}
%   \item[{\bf SM7:}] Suppose $j: A \to B$ is a cofibration and $q: X
%   \to Y$ is a fibration. Then
%   $$\xymatrix@1{{\bf Hom}_\mathcal{D}(B,X)\ar[r]^-{(j^\ast,q_\ast)}
%   &{\bf Hom}_\mathcal{D}(A,X) \times_{{\bf Hom}_\mathcal{D}(A,Y)} {\bf
%   Hom}_\mathcal{D}(B,Y)}$$
%   is a fibration of simplicial sets, which is trivial if $j$ or
%   $q$ is trivial.
%  \end{itemize}

Suppose that $\mathcal{C}$ and $\mathcal{D}$ are two closed model
categories.
%\begin{itemize}
%   \item[1.] 
We call a functor $F: \mathcal{C} \to \mathcal{D}$
a $left$ $Quillen$ $functor$ if $F$ is a left adjoint and
preserves cofibrations and trivial cofibrations.
%   \item[2.] 
We call a functor $U: \mathcal{D} \to \mathcal{C}$
   a $right$ $Quillen$ $functor$ if $U$ is a right adjoint and preserves
     fibrations and trivial fibrations.
%   \item[3.] 

Suppose that $(F,U, \varphi)$ is an adjunction from
$\mathcal{C}$ to $\mathcal{D}$. That is, $F$ is a functor
$\mathcal{C} \to \mathcal{D}$, $U$ is a functor $\mathcal{D} \to
\mathcal{C}$, and $\varphi$ is a natural isomorphism
${\mathcal{D}}(FC, D) \to {\mathcal{C}}(C, UD)$ expressing $U$ as
a right adjoint of $F$. We call $(F,U, \varphi)$ a $Quillen$
$adjunction$ if $F$ is a left Quillen functor (cf. \cite{hovey:book}).
%\end{itemize}
A Quillen adjunction $(F,U, \varphi): \mathcal{C} \to \mathcal{D}$
is called a $Quillen$ $equivalence$ if and only if, for all
cofibrant $X$ in $\mathcal{C}$ and fibrant $Y$ in $\mathcal{D}$, a
map $f: FX \to Y$ is a weak equivalence in $\mathcal{D}$ if and
only if $\varphi(f): X \to UY$ is a weak equivalence in
$\mathcal{C}$ (cf. \cite{hovey:book}).

\section{Presheaves of simplicial groupoids}

Dwyer and Kan \cite{dwyerKan:simplicialGroupoids,goerssJardine:book} show that with the following
definitions of weak equivalence, fibration and cofibration, the
category $s${\bf Gd} of simplicial groupoids satisfies the axioms
for a closed model category.

A map $f: G \to H$ of simplicial groupoids is said to be a
\emph{weak equivalence} of $s${\bf Gd} if
\begin{itemize}
 \item[(1)] the morphism $f$ induces an
isomorphism $\pi_0G \cong \pi_0H$, and
\item [(2)] each induced map
$f: G(x,x) \to H(f(x),f(x)),$ $x\in$ Ob($G$) is a weak equivalence
of simplicial groups (or of simplicial sets).
\end{itemize}
A map $g: H \to K$ of simplicial groupoids is said to be a
\emph{fibration} of $s${\bf Gd} if
\begin{itemize}
\item[(1)] the morphism $g$ has path lifting property in the sense
for every object $x$ of $H$ and morphism $\omega: g(x) \to y$ of
the groupoids $K_0$, there is a morphism $\hat{\omega}: x \to z$
of $H_0$ such that $g(\hat{\omega}) = \omega$, and
\item[(2)] each
induced map $g: H(x,x) \to K(g(x),g(x)),$ $x\in$ Ob($H$) is a
fibration of simplicial groups (or of simplicial sets).
\end{itemize}
A \emph{cofibration} of simplicial groupoids is defined to be a
map which has the left lifting property with respect to all
morphisms of $s${\bf Gd} which are both fibrations and weak
equivalences.

Let $ \mathcal{C}$ be a fixed small Grothendieck site. That is, a small category with a Grothendieck topology.
There is an adjunction between the loop groupoid functor $G: {\bf
S} \to s{\bf Gd} $ and the universal cocycle functor
$\overline{W}$ \cite[Lemma V.7.7 ]{goerssJardine:book}. By applying these
functors pointwise to simplicial presheaves and presheaves of
simplicial groupoids, one obtains functors
$$G: {\bf S}Pre(\mathcal{C})\rightleftarrows s{\bf Gd}Pre(\mathcal{C}): \overline{W}$$

So there is

%%%%%%%%%%%%%%%%%%%%%%
%  proposition 2.1   %
%%%%%%%%%%%%%%%%%%%%%%
\begin{proposition}\label{P:adjoint}
The functor $G: {\bf S}Pre(\mathcal{C}) \to s{\bf Gd}Pre(
\mathcal{C})$ is left adjoint to the functor $\overline{W}$.
\end{proposition}

A map $f: X \to Y$ in the category $s${\bf Gd}Pre($ \mathcal{C}$)
is said to be a \emph{fibration} if the induced map
$\overline{W}(f): \overline{W}X \to \overline{W}Y$ is a global
fibration in the category {\bf S}Pre($ \mathcal{C}$) in the sense
of \cite{jardine:simplicialPresheaves}.
A map $g: Z \to U$ in the category $s${\bf Gd}Pre($ \mathcal{C}$)
is said to be a \emph{weak equivalence} if the induced map
$\overline{W}(g): \overline{W}Z \to \overline{W}U$ is a
topological weak equivalence in the category {\bf S}Pre($
\mathcal{C}$) in the sense of \cite{jardine:simplicialPresheaves}.
A \emph{cofibration} in the category $s${\bf Gd}Pre($
\mathcal{C}$) is a map of presheaves of simplicial groupoids which
has the left lifting property with respect to all fibrations and
weak equivalences.
Say that a map of presheaves of simplicial groupoids $f$ is a
\emph{trivial fibration} if it is both a fibration and a weak
equivalence; a map $g$ is a \emph{trivial cofibration} if it is
both a cofibration and a weak equivalence.

Here we follow Jardine~\cite{jardine:simplicialPresheaves}. The site $\mathcal{C}$ is small, so
that there is a cardinal number $\alpha$ such that $\alpha$ is larger
than the cardinality of the set of subsets \textbf{P}Mor($
\mathcal{C}$) of the set of morphisms Mor($ \mathcal{C}$) of $
\mathcal{C}$. A simplicial presheaf $X$ is said to be $\alpha-bounded$
if the cardinality of each $X_n(U), U \in \mathcal{C}, n \geq 0$, is
smaller than $\alpha$.

A map $p: X \to Y$ in the category {\bf S}Pre($ \mathcal{C}$) is a
global fibration if and only if it has the right lifting property
with respect to all trivial cofibrations $i: U \to V$ such that
$V$ is $\alpha$-bounded \cite[Lemma 2.4]{jardine:simplicialPresheaves}. Then a map $q: G \to
H$ in the category $s${\bf Gd}Pre($ \mathcal{C}$) is a fibration
if and only if it has the right lifting property with respect to
all maps $G(i): GU \to GV$ induced by those maps $i: U \to V$
since there exist the adjoint diagrams:
$$
 \xymatrix{ GU\ar[r]\ar[d]_{G(i)} &  G\ar[d]^q &&U\ar[r]\ar[d]_i & \overline{W}G\ar[d]^{\overline{W}(q)}&&(\textbf{D})\\
          GV\ar[ur]\ar[r]    &  H &&V\ar[ur]\ar[r]     &  \overline{W}H         }
$$
For each $W \in \mathcal{C}$, $GV(W)_n$ is the free groupoid on
generators $x \in V(W)_{n+1}$ subject to some relations, and
Ob$(GV(W)) = V(W)_0$, so the cardinality of each Mor($GV(W)_n$),
$n\geq 0$ and Ob$(GV(W))$ is smaller than $\beta = max(2^\alpha,
\infty)$. We also call the presheaf of simplicial groupoids $GV$
is $\beta-bounded$.

When $G$ is a simplicial group there is a pullback diagram
$$
 \xymatrix{ G\ar[r]^-{i}\ar[d]    & WG\ar[d]^q \\
           \ast \ar[r]_-{\ast}       & \overline{W}G  }
$$
where $q$ is a fibration of simplicial sets \cite[Lemma
V.4.1]{goerssJardine:book}, $G$ is the fibre over the unique vertex $\ast \in
\overline{W}G$. $G$ is a simplicial group, so $G$ is a Kan complex
\cite[Lemma I.3.4]{goerssJardine:book}. $\overline{W}G$ is a Kan complex
\cite[Corollary V.6.8]{goerssJardine:book}, so is $WG$, then for any vertex $v
\in G$ there exists a long exact sequence
$$...\to \pi_n(G,v) \stackrel{i_*}{\longrightarrow} \pi_n(WG,v)
\stackrel{q_*}{\longrightarrow} \pi_n(\overline{W}G, \ast)
\stackrel{\partial}{\longrightarrow} \pi_{n-1}(G,v) \to ...
$$
$$...\stackrel{q_*}{\longrightarrow} \pi_1(\overline{W}G,
\ast)\stackrel{\partial}{\longrightarrow} \pi_0(G)
\stackrel{i_*}{\longrightarrow}
\pi_0(WG)\stackrel{q_*}{\longrightarrow} \pi_0(\overline{W}G)$$ by
Lemma I.7.3 in \cite{goerssJardine:book}. $WG$ is contractible \cite[Lemma
V.4.6]{goerssJardine:book}, so $\pi_n(WG,v) = 0, n \geq 1$; and $\pi_0(WG) = 0$,
since for any two vertices $a,b \in WG_0 = G_0$, there exists a
1-simplex $(s_0b, b^{-1}a) \in WG_1 = G_1 \times G_0$, s.t.,
$d_1(s_0b, b^{-1}a) = b, d_0(s_0b, b^{-1}a) = a$. Then
$$\pi_n(G,v) = \pi_{n+1}(\overline{W}G, \ast), n \geq 1$$
$$\pi_0G = \pi_1(\overline{W}G, \ast)$$

For an ordinary groupoid $H$, it's standard to write $\pi_0H$ for
the set of path components of $H$. By this one means that
$$ \pi_0 H = Ob(H)/ \sim $$
where there is a relation $x\sim y$ between two objects of $H$ if
and only if there is a morphism $x \to y$ in $H$.

If now $A$ is a simplicial groupoid, all of the simplicial
structure functors $\theta^\ast : A_n \to A_m$ induce isomorphisms
$\pi_0 A_n \cong \pi_0 A_m$. We shall therefore refer to $\pi_0
A_0$ as the set of path components of the simplicial groupoid $A$,
and denote it by $\pi_0 A$.
When $A$ is a simplicial groupoid, Ob$(A) = (\overline{W}A)_0$,
Mor$(A_0) = (\overline{W}A)_1$, so $\pi_0A \cong
\pi_0(\overline{W}A)$.
Choose a representative $x$ for each $[x]\in \pi_0A$, the
inclusion
$$i: \bigsqcup_{[x] \in \pi_0 A} A(x,x) \to A$$
is a homotopy equivalence of simplicial groupoids, and the induced
map
$$\overline{W}(i): \overline{W}(\bigsqcup_{[x] \in \pi_0 A} A(x,x)) \to \overline{W}A$$
is a weak equivalence of simplicial sets. $\overline{W}$ preserves
disjoint unions, $\overline{W}(\bigsqcup_{[x] \in \pi_0 A} A(x,x))
 = \bigsqcup_{[x] \in \pi_0 A} \overline{W}(A(x,x))$ \cite[ p. 303,304 ]{goerssJardine:book}.
$$ \pi_n(\overline{W}(\bigsqcup_{[x] \in \pi_0 A} A(x,x)), x ) =
\pi_n ( \overline{W}(A(x,x)), \ast) \cong \pi_{n-1}(A(x,x), v) , n
\geq 2, v \in A(x,x)_0 $$
$$ \pi_1(\overline{W}(\bigsqcup_{[x] \in \pi_0 A} A(x,x)), x ) =
\pi_1 ( \overline{W}(A(x,x)), \ast) \cong \pi_0(A(x,x)).$$ so one
obtains
$$\pi_n(A(x,x),v) \cong \pi_{n+1}(\overline{W}A,x), n \geq 1, x
\in \mbox{Ob}(A), v \in A(x,x)_0,$$
$$\pi_0(A(x,x)) \cong \pi_1(\overline{W}A,x).$$

According to the definition of topological weak equivalence of
simplicial presheaves in \cite{jardine:simplicialPresheaves} and the relations between
simplicial groupoids and simplicial sets, we can give an explicit
description of the weak equivalences of presheaves of simplicial
groupoids.

For any presheaf of simplicial groupoids $X$ and any object $U \in
\mathcal{C}$ and $x \in \mbox{Ob}(X(U))$, $X(U)(x,x)$ is a
simplicial group. Associated to this presheaf of simplicial
groupoids $X$ on $ \mathcal{C}$ and $\ast \in X(U)(x,x)_0$ is a
presheaf $\pi^{simp}_n(X|_U,x,\ast) (n \geq 1)$ on the comma
category $ \mathcal{C}\downarrow{U}$, the presheaf of simplicial
homotopy groups of $X|_U$, based at $\ast$, which is defined by
$$(\mathcal{C}\downarrow{U})^{op} \to {\bf Grp}$$
$$\varphi: V \to U \mapsto \pi_n(X(V)(x_V,x_V), \ast_V)$$
where $x_V$ and $\ast_V$ are the images of $x$ and $\ast$ in
$X(V)$ under the map $X(U) \to X(V)$ which is induced by $V \to
U$, respectively. The simplicial homotopy group
$\pi_n(X(V)(x_V,x_V), \ast_V)$ exists since the simplicial group
$X(V)(x_V,x_V)$ is a Kan complex \cite[Lemma I.3.4]{goerssJardine:book}.

Let $\pi_n(X|_U,x,\ast)$ be the associated sheaf of the presheaf
$\pi^{simp}_n(X|_U,x,\ast)$, i.e., $\pi_n(X|_U,x,\ast) =
L^2\pi^{simp}_n(X|_U,x,\ast)$. Then $\pi_n(X|_U,x,\ast)$ is a
sheaf of groups which is abelian if $n \geq 2$. The sheaves
$\pi_0(X|_U,x)$ and $\pi_0(X)$ of path components are defined
similarly.

A map $f: X \to Y$ of presheaves of simplicial groupoids is a
\emph{weak equivalence} if it induces isomorphisms of sheaves
$$f_*: \pi_n(X|_U,x,\ast) \cong \pi_n(Y|_U,fx,f\ast), n \geq 1, U \in \mathcal{C}, x \in \mbox{Ob}(X(U)),
\ast \in X(U)(x,x)_0$$
$$f_*: \pi_0(X|_U,x) \cong \pi_0(Y|_U,fx).$$
$$f_*: \pi_0(X) \cong \pi_0(Y).$$

In view of Proposition 1.18 in \cite{jardine:simplicialPresheaves}, the weak equivalences are
just same as the combinatorial weak equivalences in \cite{jardine:simplicialPresheaves}.  Since
the weak equivalences are given by the isomorphisms between sheaves,
thus, Proposition 1.11 in \cite{jardine:simplicialPresheaves} implies (or directly from {\bf
  CM2} for the category {\bf S}Pre($ \mathcal{C}$))

%%%%%%%%%%%%%%%%%
% lemma 3.1.1   %
%%%%%%%%%%%%%%%%%
\begin{lemma}\label{L:3we}
Given maps of presheaves of simplicial groupoids $f: X \to Y$ and
$g: Y \to Z$, if any two of $f,g$, or $g\circ f$ are weak
equivalences, then so is the third.
\end{lemma}

%%%%%%%%%%%%%%%%%
% lemma 3.1.2   %
%%%%%%%%%%%%%%%%%
\begin{lemma}\label{L:pweak}
The functor $ X \mapsto \overline{W}G(X)$ preserves weak
equivalences of simplicial presheaves.
\end{lemma}
\begin{proof}
When $T$ is a simplicial set, the natural simplicial map $\eta: T
\to \overline{W}G(T)$ is a weak equivalence of simplicial sets
\cite[Theorem V.7.8]{goerssJardine:book}. So the map $ X \to \overline{W}G(X)$ is
a pointwise weak equivalence of simplicial presheaves, then it is
a weak equivalence.

There exists a commutative diagram
$$
 \xymatrix{ X\ar[r]^-{\eta_X}\ar[d]_f    & \overline{W}G(X)
 \ar[d]^{\overline{W}G(f)}\\
 Y\ar[r]_-{\eta_Y}       & \overline{W}G(Y)  }
$$
where both $\eta_X$ and $\eta_Y$ are weak equivalences, if $f: X
\to Y$ is a weak equivalences, so is $\overline{W}G(f)$ by the
{\bf CM2} of the closed model category \textbf{S}Pre($
\mathcal{C}$).
\end{proof}

%%%%%%%%%%%%%%%%%%%
%  lemma 3.1.3    %
%%%%%%%%%%%%%%%%%%%
\begin{lemma}\label{L:pco}
The functor $G: \textbf{S}Pre(\mathcal{C}) \to
s\textbf{Gd}Pre(\mathcal{C})$ preserves cofibrations and weak
equivalences.
\end{lemma}
\begin{proof}
The adjoint diagrams (\textbf{D}) imply that the functor $G$
preserves cofibrations. Lemma~\ref{L:pweak} implies that $G$
preserves weak equivalences.
\end{proof}

%%%%%%%%%%%%%%%%%%%%%
%  lemma 3.1.4      %
%%%%%%%%%%%%%%%%%%%%%
\begin{lemma}\label{L:push}
The category $s$\textbf{Gd}Pre($\mathcal{C}$) has all pushouts,
and is hence cocomplete. The class of cofibrations in
$s$\textbf{Gd}Pre($\mathcal{C}$) is closed under pushout.
\end{lemma}
\begin{proof}
The category $s$\textbf{Gd} has all pushouts and is cocomplete, so
is the category $s$\textbf{Gd}Pre($\mathcal{C}$), since we can
take the pushout and colimit pointwise. The second statement is
obvious.
\end{proof}

There exists a Kan $Ex^\infty$ functor from {\bf
S}Pre$(\mathcal{C})$ to {\bf S}Pre$(\mathcal{C})$, such that
$Ex^\infty X$ is locally fibrant for any simplicial presheaf $X$
and the canonical map $\nu: X \to Ex^\infty X$ is a pointwise weak
equivalence \cite{jardine:booleanLocalization}.

Fix a Boolean localization $\wp:Shv(\mathcal{B}) \to
{\mathcal{E}}$, and consider the functors
$${\bf S}Pre(\mathcal{C}) \stackrel{L^2}{\longrightarrow} {\bf
S}{\mathcal{E}} \stackrel{\wp^*}{\longrightarrow} {\bf
S}Shv({\mathcal{B}}),$$ 
%relating the categories of simplicial
%presheaves on $\mathcal{C}$ and the categories of simplicial
%sheaves and the categories of simplicial objects in the categories
%of sheaves $Shv(\mathcal{B})$, 
where $L^2$ is the associated sheaf
functor. In \cite{jardine:booleanLocalization} Jardine proves that the topological weak
equivalence between simplicial presheaves \cite{jardine:simplicialPresheaves} coincides with
the local weak equivalence \cite{jardine:booleanLocalization}, i.e., a map $f: X \to Y $ of
simplicial presheaves on $\mathcal{C}$ is a topological weak
equivalence if the induced map $\wp^*L^2: \wp^*L^2Ex^\infty X \to
\wp^*L^2Ex^\infty Y$ is a pointwise weak equivalence.

Notice that there is a commutative diagram
 $$
 \xymatrix{ s{\bf Gd}Pre(\mathcal{C})\ar[r]^-{L^2}\ar[d]_{\overline{W}}& s{\bf
Gd}\mathcal{E}\ar[r]^-{\wp^*}\ar[d]^{\overline{W}}  & s{\bf
Gd}Shv(\mathcal{B}) \ar[d]^{\overline{W}}\\
        {\bf S}Pre(\mathcal{C})\ar[r]_-{L^2}& {\bf S}\mathcal{E}
         \ar[r]_-{\wp^*}  & {\bf S}Shv(\mathcal{B})   }
$$
$\overline{W}G$ is locally fibrant simplicial presheaf for any
presheaf of simplicial groupoid $G$, so a map $f: G \to H $ of
presheaves of simplicial groupoids on $\mathcal{C}$ is a weak
equivalence if the induced map $\wp^*L^2: \wp^*L^2 G \to  \wp^*L^2
H$ is a pointwise weak equivalence.

%%%%%%%%%%%%%%%%%%%%%%%
% proposition 3.1.2   %
%%%%%%%%%%%%%%%%%%%%%%%
\begin{proposition}
Trivial cofibrations of presheaves of simplicial groupoids are
closed under pushout.
\end{proposition}
\begin{proof}
Suppose that
 $$
 \xymatrix{ G\ar[r]\ar[d]_i    &  C \ar[d]^{i^\prime}\\
           H\ar[r]     &  D     }
$$
is a pushout in the category $s$\textbf{Gd}Pre($\mathcal{C}$). $i$
is a trivial cofibration, then $i^\prime$ is a cofibration by
Lemma~\ref{L:push}.

The heart of the matter for this proof is the weak equivalence.
Both $L^2$ and $\wp^*$ are left adjoint functors, so the functor
$\wp^*L^2$ preserves the pushout
 $$
 \xymatrix{ \wp^*L^2G\ar[r]\ar[d]_{\wp^*L^2(i)}    &  \wp^*L^2C
 \ar[d]^{\wp^*L^2(i^\prime)}\\
       \wp^*L^2H\ar[r]     &  \wp^*L^2D        }
$$
The map $\wp^*L^2(i)$ is a pointwise weak equivalence and
pointwise cofibration, so for any $U \in {\mathcal{B}}$, the
diagram
 $$
 \xymatrix{ \wp^*L^2G(U)\ar[r]\ar[d]_{\wp^*L^2(i)}    &  \wp^*L^2C(U)
 \ar[d]^{\wp^*L^2(i^\prime)}\\
   \wp^*L^2H(U) \ar[r]     &  \wp^*L^2D(U)  }
$$
is a pushout in the category $s${\bf Gd}. Since the category $s${\bf
  Gd} is a closed model category in which the map $\wp^*L^2(i)$ is a
trivial cofibration, the map $\wp^*L^2(i^\prime): \wp^*L^2C(U) \to
\wp^*L^2D(U)$ is also a trivial cofibration. Then $\wp^*L^2(i^\prime):
\wp^*L^2C \to \wp^*L^2D$ is a pointwise weak equivalence, and thus
$i^\prime: C \to D$ is a weak equivalence in the category $s${\bf
  Gd}Pre($\mathcal{C}$).
\end{proof}

Given a trivial cofibration $i: A \to B$ in the category
\textbf{S}Pre($\mathcal{C}$), suppose that
 $$
 \xymatrix{ GA \ar[r]\ar[d]_{G(i)}&  C \ar[d]^{i^\prime}\\
           GB  \ar[r]     &  D      }
$$
is a pushout in the category $s$\textbf{Gd}Pre($\mathcal{C}$). The
map $G(i)$ is a trivial cofibration by Lemma~\ref{L:pco}, then the
map $i^\prime$ is a trivial cofibration.

%%%%%%%%%%%%%%%%
% lemma 3.1.5  %
%%%%%%%%%%%%%%%%
\begin{lemma}\label{L:fatcf}
Every map $f: X \to Y$ of presheaves of simplicial groupoids may
be factored
$$
 \xymatrix{ X \ar[rr]^-{f}\ar[dr]_-i   & & Y \\
                 &Z\ar[ur]_-p   }
$$
where $i$ is a trivial cofibration and $p$ is a fibration.
\end{lemma}
\begin{proof}
We use the transfinite small object argument. Choose a cardinal number
$\gamma > 2^\beta$, and define a functor $F: \gamma \to s{\bf
Gd}Pre({\mathcal{C}})\downarrow Y$ on the partially ordered set
$\gamma$ by setting $F(0) = f: X \to Y, F(s): X(s) \to Y$ such
that
\begin{itemize}
  \item[(1)] $X(0) = X$,
  \item[(2)] $X(t) = \underset{ \overrightarrow{s<t}}{\lim} X(s)$ for all
            limit ordinals $t < \gamma$, and
  \item[(3)] the map $X(s) \to X(s+1)$ is defined by the pushout diagram
          $$
 \xymatrix{ \bigsqcup_D GU_D \ar[r]^-{(\alpha_D)}\ar[d]_-{\bigsqcup_D Gi_D} & X(s)
 \ar[d]\\
       \bigsqcup_D GV_D  \ar[r]       & X(s+1) }
$$
where the index $D$ refers to a set of representatives for all
diagrams
$$
 \xymatrix{  GU_D   \ar[r]^-{\alpha_D}\ar[d]_-{Gi_D} & X(s)
 \ar[d]^-{F(s)} \\
     GV_D  \ar[r]       &   Y }
$$
such that $Gi_D: GU_D \to GV_D$ is induced by $i_D: U_D \to V_D$,
where $i_D$ is a trivial cofibration in {\bf
S}Pre({$\mathcal{C}$}) with $V_D$ $\alpha-$bounded.
\end{itemize}

Then $GV_D$ is $\beta-$bounded. Let
$$ X(\gamma) = \lim_{\longrightarrow \atop {t < \gamma}} X(t)$$
and consider the induced factorization of $f$
$$
 \xymatrix{ X\ar[rr]^-{i_\gamma}\ar[dr]_-f    && X(\gamma)
 \ar[dl]^-{F(\gamma)} \\
      &    Y        }
$$
Then $i_\gamma$ is a trivial cofibration, since it is a filtered
colimit of such. Also, for any diagram
$$
 \xymatrix{ GU  \ar[d]_-{Gi}\ar[r]    & X(\gamma)
 \ar[d]^-{F(\gamma)}\\
         GV \ar[r]       & Y  }
$$
such that $GV$ is $\beta-$bounded and $Gi$ is a trivial
cofibration, the map $GU \to X(\gamma)$ must factor through some
$X(n) \to X(\gamma), n < \gamma$, for otherwise $GU$ has too many
subobjects.
\end{proof}

For each object $U$ of $\mathcal{C}$, the $U-$sections functor $X
\to X(U)$ has a left adjoint $?_U: {\bf S} \to {\bf
S}Pre(\mathcal{C})$ which sends the simplicial set $Y$ to the
simplicial presheaf $Y_U$, which is defined by
$$Y_U(V) = \coprod_{\varphi: V \to U} Y.$$
Then a simplicial presheaves map $q: Z \to X$ is a trivial
fibration if and only if it has the right lifting property with
respect to all inclusions $S\subset \triangle^n_U$ of subobjects
of the $\triangle^n_U, U \in \mathcal{C}, n \geq 0 $ \cite[p. 68
]{jardine:simplicialPresheaves}. So a map $p: G \to H$ of presheaves of simplicial groupoids
is a trivial fibration if and only if it has the right lifting
property with respect to all inclusions $GS\subset G\triangle^n_U$
of subobjects of the $G\triangle^n_U, U \in \mathcal{C}, n \geq 0
$. A transfinite small object argument, as in Lemma~\ref{L:fatcf},
shows that

%%%%%%%%%%%%%%%%%%%%%
% lemma 3.1.6       %
%%%%%%%%%%%%%%%%%%%%%
\begin{lemma}\label{L:factf}
Every map $g: Z \to W$ of presheaves of simplicial groupoids may
be factored
$$
 \xymatrix{ Z \ar[rr]^-{g}\ar[dr]_-j &    & W \\
          & M\ar[ur]_-q            }
$$
where $j$ is a cofibration and $q$ is a trivial fibration.
\end{lemma}

%%%%%%%%%%%%%%%%%%%%%%%%
%    lemma 3.1.7       %
%%%%%%%%%%%%%%%%%%%%%%%%
\begin{lemma}\label{L:lift}
For the commutative diagram
$$
 \xymatrix{ U  \ar[r]\ar[d]_-i    & X \ar[d]^-p\\
          V  \ar[ur]^-s\ar[r]       & Y  }
$$
where $i$ is a trivial cofibration and $p$ is a fibration in the
category $s${\bf Gd}Pre($\mathcal{C}$), there exists a lifting
$s$.
\end{lemma}
\begin{proof}
Suppose that $i: U \to V$ is a trivial cofibration. Then $i$ has a
factorization
$$
\xymatrix{ U\ar[r]^-j\ar[d]_-i & W \ar[dl]^-q\\
                 V            }
$$
where $q$ is a fibration and $j$ is a trivial cofibration which
has the left lifting property with respect to all fibrations by
the construction in the proof of Lemma \ref{L:fatcf}. But then $q$
is a trivial fibration, and so the lifting exists in the diagram
$$
 \xymatrix{ U\ar[r]^-j\ar[d]_-i    & W\ar[d]^-q\\
                 V  \ar[ur]\ar[r]_-{1_V}       & V  }
$$
It follows that $i$ is a retract of $j$, so that $i$ has the left
lifting property with respect to all fibrations.
\end{proof}

%%%%%%%%%%%%%%%%%%%%%%
%  theorem 3.1.1     %
%%%%%%%%%%%%%%%%%%%%%%
\begin{theorem}\label{T:main}
The category $s${\bf Gd}Pre($\mathcal{C}$), with the classes of
fibrations, weak equivalences and cofibrations as defined above,
satisfies the axioms for a closed model category.
\end{theorem}
\begin{proof}
The category $s${\bf Gd} is closed under all finite limits and
colimits. So taking the limits and colimits pointwise, the
category $s${\bf Gd}Pre($\mathcal{C}$) is also closed under all
finite limits and colimits. This is {\bf CM1}. {\bf CM2} is
Lemma~\ref{L:3we}. {\bf CM3} is trivial. The first part of {\bf
CM4} is Lemma~\ref{L:lift}, the second part is the definition
of a cofibration. {\bf CM5}(1) is Lemma~\ref{L:fatcf}, {\bf
CM5}(2) is Lemma~\ref{L:factf}.
\end{proof}

%%%%%%%%%%%%%%%%%%%%%%%%
%  remark 3.1.1        %
%%%%%%%%%%%%%%%%%%%%%%%%
\begin{remark}
Fibrations (trivial fibrations) in the category $s${\bf
Gd}Pre($\mathcal{C}$) have the right lifting property with respect
to all maps $G(i): GU \to GV$ induced by maps $i: U \to V$
where $i$ is a trivial cofibration (cofibration) in the category
{\bf S}Pre($\mathcal{C}$) and $V$ is $\alpha-$bounded. So the
category $s${\bf Gd}Pre($\mathcal{C}$) is cofibrantly generated.
\end{remark}

%%%%%%%%%%%%%%%%%%%%%%
%  lemma  3.1.8      %
%%%%%%%%%%%%%%%%%%%%%%
\begin{lemma}\label{L:preserve}
\begin{itemize}
   \item[(1)] The functor $\overline{W}: s{\bf Gd}Pre(\mathcal{C}) \to
{\bf S}Pre(\mathcal{C})$ preserves fibrations and weak
equivalences.
   \item[(2)] A map $K \to \overline{W}X \in {\bf S}Pre(\mathcal{C})$
   is a weak equivalence if and only if its adjoint $GK \to X \in
s{\bf Gd}Pre(\mathcal{C})$ is a weak equivalence.
\end{itemize}
\end{lemma}
\begin{proof}
(1) This is implied by the definitions of fibration and weak
equivalence.

(2) There is a commutative diagram
$$
 \xymatrix{ K    \ar[rr]\ar[dr]   & & \overline{W}GK\ar[dl]\\
         &         \overline{W}X         }
$$
where the map $K \to \overline{W}GK$ is a pointwise weak
equivalence \cite[Theorem V.7.8(3)]{goerssJardine:book}. So the map $K \to
\overline{W}X$ is a weak equivalence if and only if the map
$\overline{W}GK \to \overline{W}X$ is a weak equivalence, i.e.,
the map $GK \to X$ is a weak equivalence.
\end{proof}

%%%%%%%%%%%%%%%%%%%%%%%%%%%
%   theorem 3.1.2         %
%%%%%%%%%%%%%%%%%%%%%%%%%%%
\begin{corollary}\label{C:homoequ}
The functor $G$ and $\overline{W}$ induce an equivalence of
homotopy categories
$$ \mbox{Ho}(s{\bf Gd}Pre(\mathcal{C})) \simeq \mbox{Ho}({\bf S}Pre(\mathcal{C}))$$
\end{corollary}
\begin{proof}
Lemma~\ref{L:preserve} implies that the natural maps $\varepsilon:
G\overline{W}K \to K$ and $\eta: X \to \overline{W}GX$ are weak
equivalences for all presheaves of simplicial groupoids $K$ and
simplicial presheaves $X$.
\end{proof}

%%%%%%%%%%%%%%%%%%%%%%%%%%%
%   corollary 3.1.1       %
%%%%%%%%%%%%%%%%%%%%%%%%%%%
\begin{corollary}
The adjunction $G: {\bf S}Pre(\mathcal{C})\rightleftarrows s{\bf
Gd}Pre(\mathcal{C}): \overline{W}$ is a Quillen equivalence.
\end{corollary}
\begin{proof}
It follows from Theorem~\ref{T:main},
Proposition~\ref{P:adjoint}, Lemma~\ref{L:pco},
Lemma~\ref{L:preserve} and the definition of a Quillen adjunction.
\end{proof}

\begin{theorem}\label{T:simpresgpro}
The category $s${\bf Gd}Pre($ \mathcal{C}$) is right proper.
\end{theorem}
\begin{proof}
Given a pullback diagram in $s${\bf Gd}Pre($ \mathcal{C}$)
$$\xymatrix{
X \ar[r]^-{g_\ast} \ar[d] & Y\ar[d]^-p\\
Z \ar[r]_-g& W}
$$
with $p$ a fibration and $g$ a weak equivalence, there exists a
pullback diagram in {\bf S}Pre($ \mathcal{C}$)
$$\xymatrix{
\overline{W}X \ar[r]^-{\overline{W}g_\ast} \ar[d] & \overline{W}Y\ar[d]^-{\overline{W}p}\\
\overline{W}Z \ar[r]_-{\overline{W}g}& \overline{W}W}
$$
$\overline{W}$ preserves fibrations and weak equivalences, hence
$\overline{W}p$ is a fibration and $\overline{W}g$ is a weak
equivalence in {\bf S}Pre($ \mathcal{C}$). {\bf S}Pre($
\mathcal{C}$) is proper, so $\overline{W}g_\ast$ is a weak
equivalence as well, hence the map $g_\ast$ is a weak equivalence.
So the axiom {\bf P1} holds.
\end{proof}

\section{Presheaves of 2-groupoids}

% {\bf 2-Gpd}Pre($\mathcal{C}$) is the category of presheaves of
% 2-groupoids on $\mathcal{C}$; its objects are the contravariant
% functors from $\mathcal{C}$ to the category {\bf 2-Gpd} of
% 2-groupoids, and its morphisms are natural transformations.

Moerdijk and Svensson show that \cite{moerdijkSvensson:homotopy2types}, with the following
definitions of weak equivalence, fibration and cofibration, the
category {\bf 2-Gpd} of 2-groupoids satisfies the axioms for a
closed model category.

A map $\varphi: \mathcal{A} \to \mathcal{B}$ of 2-groupoids is
said to be a \emph{weak equivalence} of {\bf 2-Gpd} if
\begin{itemize}
\item[(1)] for every object $b$ of $\mathcal{B}$ there exists an
object $a$ of $\mathcal{A}$ and an arrow $\varphi(a) \to b$;
\item[(2)] for any two objects $a,a'$ in $\mathcal{A}$, $\varphi$
induces an equivalence of categories (groupoids)
$$\varphi_{a, a'}: Hom_{\mathcal{A}}(a,a') \to
Hom_{\mathcal{B}}(\varphi(a),\varphi(a')).$$
\end{itemize}

A map $\psi: \mathcal{B} \to \mathcal{A}$ of 2-groupoids is said
to be a ($Grothendieck$) \emph{fibration} of {\bf 2-Gpd} if for
any arrow $f: b_1 \to b_2$ in $\mathcal{B}$ and any arrows $g: a_0
\to \psi(b_1)$ and $h: a_0 \to \psi(b_2)$, any deformation
$\alpha: h \Rightarrow \psi(f)\circ g$ can be lifted to a
deformation $\tilde{\alpha}: \tilde{h} \Rightarrow f \circ
\tilde{g}$ in $\mathcal{B}$ (in the sense that
$\psi(\tilde{\alpha}) = \alpha, \psi(\tilde{h}) = h$ and
$\psi(\tilde{g}) = g$ ).

A \emph{cofibration} of 2-groupoids is defined to be a map which
has the left lifting property with respect to all morphisms of
{\bf 2-Gpd} which are both fibrations and weak equivalences.

There is an adjunction \cite{moerdijkSvensson:homotopy2types}:
$$G: {\bf S}\rightleftarrows {\bf 2-Gpd}: \overline{W}$$
where the functor $\overline{W}$ is the functor $N$ in \cite{moerdijkSvensson:homotopy2types}
and the functor $G$ is the Whitehead 2-groupoid functor $W$ in
\cite{moerdijkSvensson:homotopy2types}: $W(X) = W(|X|,|X^{(1)}|,|X^{(0)}|)$. By applying these
functors pointwise to simplicial presheaves and presheaves of
2-groupoids, one obtains functors
$$G: {\bf S}Pre({\bf C})\rightleftarrows {\bf 2-Gpd}Pre({\bf C}): \overline{W}$$
and there is

%%%%%%%%%%%%%%%%%%%%%
% proposition 3.2.1 %
%%%%%%%%%%%%%%%%%%%%%
\begin{proposition}\label{P:adj}
The functor $G: {\bf S}Pre({\mathcal{C}}) \to {\bf
2-Gpd}Pre(\mathcal{C})$ is left adjoint to the functor
$\overline{W}$.
\end{proposition}

A map $f: X \to Y$ in the category {\bf 2-Gpd}Pre($\mathcal{C}$)
is said to be a \emph{fibration} if the induced map
$\overline{W}(f): \overline{W}X \to \overline{W}Y$ is a global
fibration in the category {\bf S}Pre($\mathcal{C}$).
A map $g: Z \to U$ in the category {\bf 2-Gpd}Pre($\mathcal{C}$)
is said to be a \emph{weak equivalence} if the induced map
$\overline{W}(g): \overline{W}Z \to \overline{W}U$ is a weak
equivalence in the category {\bf S}Pre($\mathcal{C}$).
A \emph{cofibration} in the category {\bf 2-Gpd}Pre($\mathcal{C}$)
is a map of presheaves of 2-groupoids which has the left lifting
property with respect to all fibrations and weak equivalences.
Say that a map of presheaves of 2-groupoids $f$ is a \emph{trivial
fibration} if it is both a fibration and a weak equivalence; a map
$g$ is a \emph{trivial cofibration} if it is both a cofibration
and a weak equivalence.

%Similarly, we can define the above concepts according to the
%category $s${\bf Gd}Pre($\mathcal{C}$).
Our development follows that for presheaves of simplicial groupoids in
the previous section. We will omit proofs that are essentially the
same.

A map $q: G \to H$ in the category {\bf 2-Gpd}Pre($ \mathcal{C}$)
is a fibration if and only if it has the right lifting property
with respect to all maps $G(i): GU \to GV$ induced by the maps $i:
U \to V$ where $i$ is a trivial cofibration in the category {\bf
S}Pre($ \mathcal{C}$) and $V$ is $\alpha-$bounded, since there
exist two adjoint diagrams similar to the diagrams {\bf D}.

For each $S \in \mathcal{C}$, Ob$(GV(S)) = V(S)_0$, Mor$(GV(S))$
and 2-cell$(GV(S))$ are free generated by $V(S)_1$ and $V(S)_2$,
subject to some relations, respectively. So the cardinality of
objects, morphisms and 2-cells of 2-groupoid $GV(S)$ is smaller
than $\beta = max(2^\alpha, \infty)$, where $\alpha$ is a bound of
the simplicial presheaf $V(S)$. We also call the presheaf of
2-groupoids $GV$ is $\beta-bounded$.

For each 2-groupoid $G$ and each object $x$ of $G$, there are
natural isomorphisms \cite[Proposition 2.1(iii)]{moerdijkSvensson:homotopy2types}:
$$\pi_0(\overline{W}G) \cong \pi_0(G),$$
$$\pi_1(\overline{W}G,x) \cong \pi_1(G,x),$$
$$\pi_2(\overline{W}G,x) \cong \pi_2(G,x),$$
$$\pi_i(\overline{W}G,x) \cong 0 ~(i>2).$$

According to the definition of topological weak equivalence of
simplicial presheaves in \cite{jardine:simplicialPresheaves} and the above relations, we can
give an explicit description of weak equivalence of presheaves of
2-groupoids.

For any presheaf of 2-groupoids $X$ and any object $U \in \mathcal{C}$
and $x \in \mbox{Ob}(X(U))$, there is a presheaf on the comma category
$\mathcal{C}\downarrow{U}$, called the presheaf of homotopy groups of
$X|_U$, based at $x$, which is defined by
$$(\mathcal{C}\downarrow{U})^{op} \to {\bf Grp}$$
$$\varphi: V \to U \mapsto \pi_i(X(V),x_V), i = 1,2$$
where $x_V$ is the image of $x$ in $X(V)$ under the map $X(U) \to
X(V)$ which is induced by $V \to U$.

Let $\pi_i(X|_U,x), i = 1,2$ be the associated sheaves of the
above presheaves. The sheaf $\pi_0(X)$ of path components is
defined similarly.

A map $f: X \to Y$ of presheaves of 2-groupoids is a
\emph{weak equivalence} if it induces isomorphisms of sheaves
$$f_*: \pi_i(X|_U,x) \cong \pi_i(Y|_U,fx), i = 1,2;
U \in \mathcal{C}, x \in \mbox{Ob}(X(U))$$
$$f_*: \pi_0(X) \cong \pi_0(Y).$$

In parallel with the corresponding arguments for presheaves of
simplicial groupoids, we have

%%%%%%%%%%%%%%%%%
% lemma  3.2.1  %
%%%%%%%%%%%%%%%%%
\begin{lemma}
Given maps of presheaves of 2-groupoids $f: X \to Y$ and $g: Y \to
Z$, if any two of $f,g$, or $g\circ f$ are weak equivalences, then
so is the third.
\end{lemma}

%%%%%%%%%%%%%%%%%%
%  lemma 3.2.2   %
%%%%%%%%%%%%%%%%%%
\begin{lemma}
The functor $ X \mapsto \overline{W}G(X)$ preserves weak
equivalences of simplicial presheaves.
\end{lemma}
\begin{proof}
When $T$ is a simplicial set, there are isomorphisms \cite{moerdijkSvensson:homotopy2types}
$$\pi_0(\overline{W}GT) \cong \pi_0(GT) \cong \pi_0(T),$$
$$\pi_i(\overline{W}GT, t_0) \cong \pi_i(GT, t_0) \cong \pi_i(T, t_0)~(i=1,2),t_0 \in T_0.$$
$$\pi_i(\overline{W}GT, t_0) = 0 ~(i > 2).$$
so there exist isomorphisms of sheaves
$$\pi_0(\overline{W}GX) \cong \pi_0(GX) \cong \pi_0(X),$$
$$\pi_i(\overline{W}GX|_U,x) \cong \pi_i(GX|_U, x) \cong \pi_i(X|_U, x)~(i=1,2)~U \in  \mathcal{C}~x \in X(U)_0.$$
and $\pi_i(\overline{W}GX|_U, x) = 0~(i>2)$.
\end{proof}

%%%%%%%%%%%%%%%%%%%
%  lemma 3.2.3    %
%%%%%%%%%%%%%%%%%%%
\begin{lemma}\label{L:Gpre}
The functor $G: \textbf{S}Pre(\mathcal{C}) \to
\textbf{2-Gpd}Pre(\mathcal{C})$ preserves cofibrations and weak
equivalences.
\end{lemma}

%%%%%%%%%%%%%%%%%%
%  lemma 3.2.4   %
%%%%%%%%%%%%%%%%%%
\begin{lemma}
The category \textbf{2-Gpd}Pre($\mathcal{C}$) has all pushouts,
and is hence cocomplete. The class of cofibrations in
\textbf{2-Gpd}Pre($\mathcal{C}$) is closed under pushout.
\end{lemma}

Notice that there is a commutative diagram
 $$
 \xymatrix{ {\bf 2-Gpd}Pre(\mathcal{C}) \ar[r]^-{L^2}\ar[d]_-{\overline{W}}& {\bf
2-Gpd}\mathcal{E} \ar[r]^-{\wp^*}\ar[d]^-{\overline{W}}  & {\bf
2-Gpd}Shv(\mathcal{B}) \ar[d]^-{\overline{W}}\\
{\bf S}Pre(\mathcal{C})\ar[r]_-{L^2}& {\bf S}\mathcal{E}
         \ar[r]_-{\wp^*}  & {\bf S}Shv(\mathcal{B})       }
$$
$\overline{W}G$ is a locally fibrant simplicial presheaf for any
presheaf of 2-groupoids $G$, so a map $f: G \to H $ of presheaves
of 2-groupoids on $\mathcal{C}$ is a weak equivalence if the
induced map $\wp^*L^2: \wp^*L^2 G \to  \wp^*L^2 H$ is a pointwise
weak equivalence.

%%%%%%%%%%%%%%%%%%%%%%%
%  proposition 3.2.2  %
%%%%%%%%%%%%%%%%%%%%%%%
\begin{proposition}
Trivial cofibrations of presheaves of 2-groupoids are closed under
pushout.
\end{proposition}

Given a trivial cofibration $i: A \to B$ in the category
\textbf{S}Pre($\mathcal{C}$), suppose that
 $$
 \xymatrix{ GA   \ar[r]\ar[d]_-{G(i)}    &  C \ar[d]^-{i^\prime}\\
          GB   \ar[r]     &  D   }
$$
is a pushout in the category \textbf{2-Gpd}Pre($\mathcal{C}$).
Then the map $i^\prime$ is a trivial cofibration.

%%%%%%%%%%%%%%%%%%%%%
%  lemma 3.2.5      %
%%%%%%%%%%%%%%%%%%%%%
\begin{lemma}\label{L:fatcf2}
Every map $f: X \to Y$ of presheaves of 2-groupoids may be
factored
$$
 \xymatrix{ X  \ar[rr]^-{f}\ar[dr]_-i  &  & Y \\
             &  Z  \ar[ur]_-p         }
$$
where $i$ is a trivial cofibration and $p$ is a fibration.
\end{lemma}

A map $p: G \to H$ of presheaves of 2-groupoids is a trivial
fibration if and only if it has the right lifting property with
respect to all inclusions $GS\subset G\triangle^n_U$ of subobjects
of the $G\triangle^n_U, U \in \mathcal{C}, n \geq 0 $. A
transfinite small object argument, as in Lemma~\ref{L:fatcf},
shows that

%%%%%%%%%%%%%%%%%%%%%%%%%
%   lemma 3.2.6         %
%%%%%%%%%%%%%%%%%%%%%%%%%
\begin{lemma}\label{L:factf2}
Every map $g: Z \to W$ of presheaves of 2-groupoids may be
factored
$$
 \xymatrix{ Z \ar[rr]^-{g}\ar[dr]_-j &   & W \\
       &         M \ar[ur]_-q       }
$$
where $j$ is a cofibration and $q$ is a trivial fibration.
\end{lemma}

%%%%%%%%%%%%%%%%%%%%%%%%%
%  lemma 3.2.7          %
%%%%%%%%%%%%%%%%%%%%%%%%%
\begin{lemma}\label{L:lift2}
For the commutative diagram
$$
 \xymatrix{ U    \ar[r]\ar[d]_-i    & X \ar[d]^-p\\
              V   \ar[ur]^-s \ar[r]       & Y  }
$$
where $i$ is a trivial cofibration and $p$ is a fibration in the
category {\bf 2-Gpd}Pre($\mathcal{C}$), there exists a lifting
$s$.
\end{lemma}

%%%%%%%%%%%%%%%%%%%%%%
% theorem 3.2.1      %
%%%%%%%%%%%%%%%%%%%%%%
\begin{theorem}\label{T:main2}
The category {\bf 2-Gpd}Pre($\mathcal{C})$, with the classes of
fibrations, weak equivalences and cofibrations as defined above,
satisfies the axioms for a closed model category.
\end{theorem}

%%%%%%%%%%%%%%%%%%%%%%%%%%%
%   lemma 3.2.8           %
%%%%%%%%%%%%%%%%%%%%%%%%%%%
\begin{lemma}
\begin{itemize}
   \item[(1)] The functor $\overline{W}: {\bf 2-Gpd}Pre(\mathcal{C}) \to {\bf
S}Pre(\mathcal{C})$ preserves fibrations and weak equivalences.
   \item[(2)] The functors $G$ and $\overline{W}$ induce adjoint functors
  $$ G: Ho({\bf S}Pre(\mathcal{C})) \rightleftarrows Ho({\bf 2-Gpd}Pre(\mathcal{C})) : \overline{W}$$
  at the level of homotopy categories.
\end{itemize}
\end{lemma}
\begin{proof}
(1) This is implied by the definitions of fibration and weak
equivalence.

(2) The functors $\overline{W}$ and $G$ both preserve weak
equivalences ((1) of this Lemma and Lemma~\ref{L:Gpre}), they
localize to functors of homotopy categories. The triangular
identities for the unit and counit will still hold after
localization.
\end{proof}

%%%%%%%%%%%%%%%%%%%%%%%%%%%%%
%   corollary 3.2.1         %
%%%%%%%%%%%%%%%%%%%%%%%%%%%%%
\begin{corollary}
The adjunction $ G: {\bf S}Pre(\mathcal{C}) \rightleftarrows {\bf
2-Gpd}Pre(\mathcal{C}) : \overline{W}$ is a Quillen adjunction.
\end{corollary}
\begin{proof}
This follows from Theorem~\ref{T:main2}, Proposition~\ref{P:adj}, and
Lemma~\ref{L:Gpre}. % and the definition of Quillen adjunction.
\end{proof}

Define the category $2-types{\bf S}Pre(\mathcal{C})$ of homotopy
2-types to be the full subcategory of Ho({\bf
S}Pre($\mathcal{C}$)) given by those simplicial presheaves with
sheaves $\pi_i(X|_U, x) = 0$ for any integer $i > 2$, any object
$U \in \mathcal{C}$ and any basepoint $x\in X(U)_0$.

%%%%%%%%%%%%%%%%%%%%
%  theorem 3.2.2   %
%%%%%%%%%%%%%%%%%%%%
\begin{theorem}\label{T:ho2-type}
The functors $G$ and $\overline{W}$ induce an equivalence of
homotopy categories
$$ \mbox{Ho}({\bf 2-Gpd}Pre(\mathcal{C})) \simeq \mbox{2-types}{\bf S}Pre(\mathcal{C})$$
\end{theorem}
\begin{proof}
  For a simplicial presheaf $X$, the natural map $\eta: X \to
  \overline{W}G(X)$ is a weak equivalence if and only if
  $\pi_i(X|_U,x) = 0$ for all $i > 2$, $U \in \mathcal{C}$, and $x \in
  X(U)_0$. For any presheaf of 2-groupoids $K$,
  $\pi_i(\overline{W}K|_U, \ast) = 0$, for all $i > 2$, $U \in
  \mathcal{C}$, and $\ast \in Ob(K(U))$, and the natural map $\varphi:
  G\overline{W}(K) \to K$ is a weak equivalence.
\end{proof}

% \bibliographystyle{alpha}
% \bibliography{my}

\begin{thebibliography}{Jar06b}

\bibitem[Cra95]{crans:sheaves}
Sjoerd~E. Crans.
\newblock Quillen closed model structures for sheaves.
\newblock {\em J. Pure Appl. Algebra}, 101(1):35--57, 1995.

\bibitem[DK84]{dwyerKan:simplicialGroupoids}
W.~G. Dwyer and D.~M. Kan.
\newblock Homotopy theory and simplicial groupoids.
\newblock {\em Nederl. Akad. Wetensch. Indag. Math.}, 46(4):379--385, 1984.

\bibitem[GJ99]{goerssJardine:book}
Paul~G. Goerss and John~F. Jardine.
\newblock {\em Simplicial homotopy theory}, volume 174 of {\em Progress in
  Mathematics}.
\newblock Birkh\"auser Verlag, Basel, 1999.

\bibitem[Hov99]{hovey:book}
Mark Hovey.
\newblock {\em Model categories}, volume~63 of {\em Mathematical Surveys and
  Monographs}.
\newblock American Mathematical Society, Providence, RI, 1999.

\bibitem[Jar87]{jardine:simplicialPresheaves}
J.~F. Jardine.
\newblock Simplicial presheaves.
\newblock {\em J. Pure Appl. Algebra}, 47(1):35--87, 1987.

\bibitem[Jar96]{jardine:booleanLocalization}
J.~F. Jardine.
\newblock Boolean localization, in practice.
\newblock {\em Doc. Math.}, 1:No.\ 13, 245--275 (electronic), 1996.

\bibitem[Jar06a]{jardine:cocycle}
J.F. Jardine.
\newblock Cocycle categories.
\newblock arXiv:math.AT/0605198, 2006.

\bibitem[Jar06b]{jardine:gerbes}
J.F. Jardine.
\newblock Homotopy classification of gerbes.
\newblock arXiv:math.AT/0605200, 2006.

\bibitem[JL06]{jardineLuo}
J.~F. Jardine and Z.~Luo.
\newblock Higher principal bundles.
\newblock {\em Math. Proc. Cambridge Philos. Soc.}, 140(2):221--243, 2006.

\bibitem[Joy84]{joyal:letter}
Andr\'e Joyal.
\newblock Homotopy theory of simplicial sheaves.
\newblock unpublished (circulated as a letter to Grothendieck), 1984.

\bibitem[JT91]{joyalTierney:strongStacks}
Andr{\'e} Joyal and Myles Tierney.
\newblock Strong stacks and classifying spaces.
\newblock In {\em Category theory ({C}omo, 1990)}, volume 1488 of {\em Lecture
  Notes in Math.}, pages 213--236. Springer, Berlin, 1991.

\bibitem[JT96]{joyalTierney:simplicialGroupoids}
Andr{\'e} Joyal and Myles Tierney.
\newblock On the homotopy theory of sheaves of simplicial groupoids.
\newblock {\em Math. Proc. Cambridge Philos. Soc.}, 120(2):263--290, 1996.

\bibitem[MS93]{moerdijkSvensson:homotopy2types}
Ieke Moerdijk and Jan-Alve Svensson.
\newblock Algebraic classification of equivariant homotopy {$2$}-types. {I}.
\newblock {\em J. Pure Appl. Algebra}, 89(1-2):187--216, 1993.

\bibitem[Qui67]{quillen:homotopicalAlgebra}
Daniel~G. Quillen.
\newblock {\em Homotopical algebra}.
\newblock Lecture Notes in Mathematics, No. 43. Springer-Verlag, Berlin, 1967.

\bibitem[Qui69]{quillen:rht}
Daniel Quillen.
\newblock Rational homotopy theory.
\newblock {\em Ann. of Math. (2)}, 90:205--295, 1969.

\end{thebibliography}

\end{document}